\documentclass{article}

\usepackage{graphicx}        
\usepackage{multicol}        
\usepackage[bottom]{footmisc}

\usepackage{times}
\usepackage{mathrsfs}
\usepackage{amsmath}
\usepackage{amssymb}

\title{An Optimal Control Derivation of Nonlinear Smoothing Equations}
\author{Jin Won Kim$^1$
\and Prashant G. Mehta$^1$}
\date{}

\def\Re{\mathbb{R}}

\def\notes#1{\marginpar{\tiny #1}\typeout{Notes!
Notes!
Notes!
}}
\renewcommand{\notes}[1]{\typeout{notes!}}

\def\half{{\mathchoice{\FRAC{1}{1}{2}}%
{\FRAC{2}{1}{2}}%
{\FRAC{3}{1}{2}}%
{\FRAC{4}{1}{2}}}}

\newcommand{\tr}{\mbox{tr}}

\def\Re{\field{R}}

\def\clZ{{\cal Z}}

\def\E{{\sf E}}

\def\clZ{{\cal Z}}

\def\beq{\begin{eqnarray}} 
\def\bc{\begin{center}} 
\def\be{\begin{enumerate}}
\def\bi{\begin{itemize}} 
\def\bs{\begin{small}}
\def\bS{\begin{slide}}
\def\ec{\end{center}} 
\def\ee{\end{enumerate}}
\def\ei{\end{itemize}}
\def\es{\end{small}}
\def\eS{\end{slide}}
\def\eeq{\end{eqnarray}}

\newcommand{\ud}{\,\mathrm{d}}

\def\Re{\mathbb{R}}
\def\E{{\sf E}}

\def\clZ{{\cal Z}}

\renewcommand{\Re}{\mathbb{R}}

\def\dv{\operatorname{diag}}

\def\sJ{{\sf J}}

\def\normfactor{C}
\def\divg{\operatorname{div}}

\renewcommand{\half}{{\frac{1}{2}}}

\newtheorem{remark}{Remark}
\newtheorem{proposition}{Proposition}

\begin{document}

\maketitle

\footnotetext[1]{Jin Won Kim and Prashant G.~Mehta are with Coordinated Science Laboratory and Department of Mechanical Science and Engineering, University of Illinois at Urbana-Champaign, Urbana, IL 61801 USA 
}

\begin{abstract}The purpose of this paper is to review and highlight some connections
	between the problem of nonlinear smoothing and optimal control of the
	Liouville equation.  
	The latter has been an active area of recent research interest owing
	to work in mean-field games and optimal transportation theory.	
	The nonlinear smoothing problem is considered here for continuous-time
	Markov processes.  The observation process is modeled as a nonlinear
	function of a hidden state with an additive Gaussian measurement
	noise.  A variational formulation is described based upon the relative
	entropy formula introduced by Newton and Mitter~\cite{mitter2003}.  The resulting
	optimal control problem is formulated on the space of
	probability distributions.  The Hamilton's equation of the optimal
	control are related to the Zakai equation of nonlinear smoothing via
	the log transformation.  The overall procedure is shown to generalize the
	classical Mortensen's minimum energy estimator for the linear Gaussian
	problem.    
\end{abstract}

\vspace{0.5in}

\noindent 
\large \emph{To Michael Dellnitz on the occasion of his 60th birthday}




\section{Introduction}
\label{sec:intro}



There is a fundamental dual relationship between estimation and
control.  The most basic of these relationships is the well
known duality between controllability and observability of a linear
system~\cite[Ch. 15]{Kailath:233814}.  The relationship suggests that
the problem of filter (estimator) design can be re-formulated as a
variational problem of
optimal control.  Such variational formulations are referred to as
the duality principle of optimal filtering.  The first duality
principle appears in the seminal (1961) paper of Kalman-Bucy, where
the problem of minimum variance estimation is shown to be dual to a
linear quadratic optimal control problem.  In these classical 
settings, the dual variational formulations are of the following two types
\cite[Sec. 7.3]{bensoussan2018estimation}: (i) minimum variance
estimator and (ii) minimum energy estimator.

The classical minimum energy estimator represents a solution of the smoothing
problem.  The estimator is modeled as a controlled version
of the state process in which the process noise term is replaced by a
control input.   The optimal control input is obtained by maximizing
the log of the conditional (smoothed) distribution.  For this reason,
the estimator is also referred to as the  
maximum a posteriori (MAP) estimator.
The MAP solution coincides with the optimal smoother in the linear-Gaussian case.  The earliest construction of the minimum energy estimator is
due to Mortensen~\cite{mortensen1968}.    

A variational formulation of the nonlinear smoothing problem -- the focus of this paper -- leading
to the conditional distribution appears
in~\cite{mitter2003}.  The formulation is based upon the
variational Kallianpur-Striebel formula~\cite[Lemma
2.2.1]{van2006filtering}.
The divergence is expressed as an optimal control objective function
which turns out to be identical to the objective function considered in
the MAP estimator~\cite{mortensen1968}. 
The difference is that the constraint now
is a controlled stochastic process, in contrast to a single trajectory
in the MAP estimator.  With the optimal control input, the law of the
stochastic process is the conditional distribution.  

The purpose of this paper is to review and highlight some connections
between nonlinear smoothing and optimal control problems involving
control of probability densities.  In recent years, there has been a
lot of interest in mean-field-type optimal control problems
where the constraint is a controlled Liouville or a Fokker-Plank
equation describing the evolution of the probability
density~\cite{bensoussan2013mean,brockett2007optimal,carmona2018probabilistic}.
In this paper, it is shown that the variational formulation proposed
in~\cite{mitter2003} is easily described and solved in these terms.  The formulation as a mean-field-type optimal control problem
is more natural compared to a stochastic optimal control
formulation considered in~\cite{mitter2003}.  
In particular, the solution with the density constraint directly 
leads to the forward-backward equation of
pathwise smoothing.  This also makes explicit the connection to the log
transformation which is known to transform the Bellman equation of optimal control into the Zakai
equation of filtering~\cite{FlemingMitter82,kappen2016adaptive}.  Apart from the case of the It\^o-diffusion, the
continuous-time Markov chain is also described.  The overall procedure
is shown to generalize the classical Mortensen's minimum energy
estimator for the linear Gaussian problem.

The outline of the remainder of this chapter is as follows: the
smoothing problem and its solution in terms of the forward-backward Zakai
equation and their pathwise representation is reviewed
in~Sec.~\ref{sec:prelim}.  The variational formulation leading to a
mean-field optimal control problem and its solution appears
in~Sec.~\ref{sec:control-problem}.  The relationship to the log
transformation and to the minimum energy estimator is described.  The
conclusions appear in Sec.~\ref{sec:conc}. 
All the proofs are contained in the Appendix.

\paragraph{Notation} We denote the $i^{\text{th}}$ element of a vector by
$[\,\cdot\,]_i$, and similarly, $(i,j)$ element of a matrix is denoted
by $[\,\cdot\,]_{ij}$. $C^k(\Re^d\,;S)$ is the space of functions with
continuous $k$-th order derivative. For a function $f\in
C^2(\Re^d\,;\Re)$, $\nabla f$ is the gradient vector and $D^2 f$ is
the Hessian matrix. For a vector field $F\in C^1(\Re^d\,;\Re^d)$,
$\divg(F)$ denotes the divergence of $F$. For a vector $v\in\Re^d$,
$\dv(v)$ denotes a diagonal matrix with diagonal entries given by the
vector; $e^v$ and $v^2$ are defined in an element-wise manner, that is,
$[e^v]_i = e^{[v]_i}$ and $[v^2]_i = ([v]_i)^2$ for $i =
1,\ldots,d$. For a matrix, $\tr(\cdot)$ denotes the trace.

\section{Preliminaries and Background}
\label{sec:prelim}

\subsection{The smoothing problem}

Consider a pair of continuous-time stochastic processes $(X,Z)$.  The
state $X=\{X_t:t\in[0,T]\}$ is a Markov process taking values in the
state space $\mathbb{S}$. The observation process $Z =
\{Z_t:t\in[0,T]\}$ is defined according to the model:
\begin{equation}\label{eq:obs_model}
Z_t = \int_0^t h(X_s) \ud s + W_t
\end{equation}
where $h:\mathbb{S}\rightarrow \mathbb{R}$ is the observation function
and $W=\{W_t: t\ge 0\}$ is a standard Wiener process.

The smoothing problem is to compute the posterior distribution ${\sf P}(X_t\in \;\cdot\;|\clZ_T)$ for arbitrary $t\in [0,T]$, where
$\clZ_T := \sigma(Z_s : 0\le s \le T)$ is the sigma-field generated by
the observation up to the terminal time $T$.

\subsection{Solution of the smoothing problem}

The smoothing problem requires a model of the Markov process $X$.  In
applications involving nonlinear smoothing, a common model is
the It\^o-diffusion in Euclidean settings:


\paragraph{Euclidean state space} The state space $\mathbb{S}=\Re^d$. The
state process $X$ is modeled as an It\^o diffusion:
\begin{equation*}\label{eq:dyn_sde}
\ud X_t = a (X_t) \ud t + \sigma(X_t) \ud B_t, \quad X_0\sim \nu_0
\end{equation*}
where $a\in C^1(\Re^d; \Re^d)$, $\sigma\in C^2(\Re^d; \Re^{d\times
	p})$ and $B=\{B_t:t\ge 0\}$ is a standard Wiener process.  The initial
distribution of $X_0$ is denoted as $\nu_0(x) \ud x$ where $\nu_0(x)$
is the probability density with respect to the Lebesgue measure.
For~\eqref{eq:obs_model}, the observation function $h\in C^2(\Re^d;\Re)$.  
It is assumed that $X_0, B, W$ are mutually independent.

The infinitesimal generator of $X$, denoted as ${\cal A}$,
acts on $C^2$ functions in its domain according to 
\[
({\cal A} f)(x):= a^\top(x) \nabla f(x) + \half \tr\big(\sigma\sigma^\top(x)(D^2f)(x)\big).
\]

The adjoint operator is denoted by ${\cal A}^\dagger$.  It acts on
$C^2$ functions in its domain according to 
\[
({\cal A}^\dagger f)(x) = -\divg(af)(x) + \half \sum_{i,j=1}^d \frac{\partial^2}{\partial x_i \partial x_j}\big([\sigma\sigma^\top]_{ij} f\big)(x) .
\]


The solution of the smoothing problem is described by a
forward-backward system of stochastic partial differential equations
(SPDE) (see~\cite[Thm. 3.8]{pardoux1981non}):
\begin{subequations}\label{eq:Zakai}
	\begin{flalign}
	&\text{(forward)}:&\ud p_t (x) &= ({\cal A}^\dagger p_t)(x) \ud t +
	h(x) p_t(x) \ud Z_t & \nonumber \\
	& & p_0(x) &= \nu_0(x), \quad \forall x\in\Re^d & \label{eq:Zakai-a}\\
	&\text{(backward)}:&-\ud q_t(x) &= {(\cal A} q_t)(x) \ud t + h(x)
	q_t(x) \overleftarrow{\ud Z}_t &\nonumber \\
	& &q_T (x) &\equiv 1 &\label{eq:Zakai-b}
	\end{flalign}
\end{subequations}
where $\overleftarrow{\ud Z}_t$ denotes a backward It\^{o} integral (see \cite[Remark 3.3]{pardoux1981non}). The smoothed distribution is then obtained as follows:
\[
{\sf P}(X_t \in \ud x\; |\clZ_T) \propto p_t(x)q_t(x) \ud x.
\]
Each of~\eqref{eq:Zakai} is referred to as the Zakai equation of nonlinear
filtering.  

\subsection{Path-wise representation of the Zakai equations}

There is a representation of the forward-backward SPDEs where the only
appearance of randomness is in the coefficients.  This is referred to
as the pathwise (or robust) form of the
filter~\cite[Sec. VI.11]{rogers2000diffusions}.

Using It\^o's formula for $\log p_t$,
\begin{align*}
\ud (\log p_t) (x)&=\frac{1}{p_t(x)} ({\cal A}^\dagger p_t)(x) \ud t + h(x)\ud Z_t - \half h^2(x) \ud t.
\end{align*}
Therefore, upon defining $\mu_t(x) := \log p_t(x) - h(x)Z_t$, the forward
Zakai equation~\eqref{eq:Zakai-a} is transformed into a parabolic partial differential
equation (pde):
\begin{align}
\frac{\partial \mu_t}{\partial t} (x) &=
e^{-(\mu_t(x)+Z_th(x))}\big({\cal A}^\dagger
e^{(\mu_t(\cdot) +Z_th(\cdot) )}\big)(x) - \half h^2(x) \nonumber\\
\mu_0(x) &= \log \nu_0(x),\quad \forall x\in\Re^d . \label{eq:pw-Zakai-a}
\end{align}
Similarly, upon defining $\lambda_t(x) = \log q_t(x) + h(x) Z_t$, the
backward Zakai equation~\eqref{eq:Zakai-b} is transformed into the
parabolic pde:
\begin{align}
-\frac{\partial \lambda_t}{\partial t} (x) &= e^{-(\lambda_t(x) - Z_th(x))}\big({\cal A}e^{\lambda_t(\cdot) - Z_th(\cdot)}\big)(x) - \half h^2(x) \nonumber\\
\lambda_T(x) &= Z_Th(x),\quad \forall x\in\Re^d. \label{eq:pw-Zakai-b}
\end{align}
The pde~\eqref{eq:pw-Zakai-a}-\eqref{eq:pw-Zakai-b} are referred to as
pathwise equations of nonlinear smoothing.  

\subsection{The finite state-space case}

Apart from It\^{o}-diffusion, another common model is a Markov chain in
finite state-space settings: 

\paragraph{Finite state space} Let the state-space be $\mathbb{S} =
\{e_1,e_2,\hdots,e_d\}$, the canonical basis in $\Re^d$.  
For~\eqref{eq:obs_model}, the linear observation model is chosen
without loss of generality: for any function $h:\mathbb{S}\to\Re$, we
have $h(x)=\tilde{h}^\top x$ where $\tilde{h}\in\Re^{d}$ is defined by
$\tilde{h}_i = h(e_i)$.  
Thus, the function space on $\mathbb{S}$ is
identified with $\Re^d$.  With a slight abuse of notation, we will
drop the tilde and simply write $h(x)=h^\top x$. 

The state process $X$ is a continuous-time Markov chain evolving in $\mathbb{S}$. 
The initial distribution for
$X_0$ is denoted as $\nu_0$.  It is an element of the probability simplex
in $\Re^d$. 
The generator of the chain is denoted as $A$.  It is a
$d\times d$ row-stochastic matrix.  It acts on a function $f\in\Re^d$
through right multiplication: $f\mapsto Af$.  The adjoint operator is the
matrix transpose $A^\top$.  
It is assumed that $X$ and $W$ are mutually independent.

The solution of the smoothing problem for the finite state-space
settings is entirely analogous: Simply replace the generator ${\cal
	A}$ in~\eqref{eq:Zakai} by the matrix $A$, and the probability density by the probability mass function.  The Zakai pde is now the
Zakai sde.  The formula for the pathwise representation are also
entirely analogous:
\begin{align}
\Big[\frac{\ud \mu_t}{\ud t}\Big]_i &= [e^{-(\mu_t+Z_th)}]_i[A^\top e^{\mu_t+Z_th}]_i - \half [h^2]_i \label{eq:pw-Zakai-finite-a}\\
-\Big[\frac{\ud \lambda_t}{\ud t}\Big]_i &= [e^{-(\lambda_t-Z_th)}]_i[Ae^{\lambda_t-Z_th}]_i - \half [h^2]_i \label{eq:pw-Zakai-finite-b}
\end{align}
with boundary condition $[\mu_0]_i = \log [\nu_0]_i$ and $[\lambda_0]_i = Z_T[h]_i$, for $i=1, \ldots,d$.

\section{Optimal Control Problem}\label{sec:control-problem}

\subsection{Variational formulation}

For the smoothing problem, an optimal control formulation is derived in the following two steps:

\paragraph{Step 1} A control-modified version of the Markov process $X$ is introduced.  The controlled process is denoted as $\tilde{X} :=
\{\tilde{X}_t:0\le t\le T\}$.  
The control problem is to pick (i) the initial distribution
$\pi_0\in{\cal P}(\mathbb{S})$ and (ii) the state transition, such that
the distribution of $\tilde{X}$ equals the conditional distribution.
For this purpose, an optimization problem is formulated in the next step.

\paragraph{Step 2} The optimization problem is formulated on the space of
probability laws. 
Let $P$ denote the law for $X$, $\tilde{P}$ denote the law for
$\tilde{X}$, and $Q^z$ denote the law for $X$ given an observation path
$z=\{z_t:0\le t\le T\}$.  Assuming
these are equivalent, the objective function is the relative entropy between $\tilde{P}$ and $Q^z$:
\begin{equation*}\label{eq:rel-entropy-cost}
\min_{\tilde{P}} \quad \E_{\tilde{P}}\Big(\log \frac{\ud \tilde{P}}{\ud P}\Big) - \E_{\tilde{P}}\Big(\log\frac{\ud Q^z}{\ud P}\Big).
\end{equation*}

Upon using the Kallianpur-Striebel formula
(see~\cite[Lemma 1.1.5 and Prop. 1.4.2]{van2006filtering}), the
optimization problem is equivalently expressed as follows:
\begin{equation}\label{eq:cost-equiv-form}
\min_{\tilde{P}} \quad {\sf D}(\tilde{P}\|P) +\E\Big(\int_0^T z_t \ud h(\tilde{X}_t) + \half|h(\tilde{X}_t)|^2 \ud t- z_Th(\tilde{X}_T)\Big).
\end{equation}
The first of these terms depends upon the details of the model used to
parametrize the controlled Markov process $\tilde{X}$.  For the two
types of Markov processes, this is discussed in the following sections.    

\begin{remark}
	The Schr\"{o}dinger bridge problem is a closely related problem of
	recent research interest where one picks $\tilde{P}$ to minimize ${\sf D}(\tilde{P}\|P)$
	subject to the constraints on marginals at time $t=0$ and $T$; cf.,~\cite{chen2016relation} where connections to
	stochastic optimal control theory are also described.  Applications of
	such models to the filtering and smoothing problems is discussed
	in~\cite{reich2019data}.  There are two differences between the
	Schr\"{o}dinger bridge problem and the smoothing problem considered
	here:
	\begin{enumerate}
		\item The objective function for the smoothing problem also includes
		an additional integral term in~\eqref{eq:cost-equiv-form} to account for
		conditioning due to observations $z$ made over time $t\in[0,T]$;
		\item The constraints on the marginals at time $t=0$ and $t=T$ are not
		present in the smoothing problem.
		Rather, one is allowed to pick the initial distribution $\pi_0$ for
		the controlled process and there is no constraint present on the 
		distribution at the terminal time $t=T$.  
	\end{enumerate}
	%
\end{remark}

\subsection{Optimal control: Euclidean state-space}
\label{sec:euclidean_ss}

The modified process $\tilde{X}$ evolves on
the state space $\Re^d$.  It is modeled as a controlled
It\^{o}-diffusion
\[
\ud \tilde{X}_t = a(\tilde{X}_t) \ud t + \sigma(\tilde{X}_t)\big(u_t(\tilde{X}_t)\ud t + \ud \tilde{B}_t\big),\quad \tilde{X}_0\sim \pi_0
\]
where $\tilde{B} = \{\tilde{B}_t:0\leq t\leq T\}$ is a copy of the process noise $B$. The controlled process is parametrized by:
\begin{enumerate}
	\item The initial density $\pi_0(x)$.
	\item The control function $u\in C^1([0,T]\times\Re^d;\Re^p)$.  The
	function of two arguments is denoted as $u_t(x)$. 
\end{enumerate}
The parameter $\pi_0$ and the function $u$ are chosen as a solution of an optimal control problem.

For a given function $v\in C^1(\Re^d;\Re^p)$, the generator of the controlled
Markov process is denoted by $\tilde{\cal A}(v)$.  It acts on a $C^2$
function $f$ in its domain according to
\begin{align*}
(\tilde{\cal A}(v)f)(x)
&= ({\cal A}f)(x) + (\sigma v)^\top(x)\nabla f(x).
\end{align*}
The adjoint operator is denoted by ${\cal A}^\dagger(v)$.  It acts on
$C^2$ functions in its domain according to 
\begin{align*}
(\tilde{\cal A}^\dagger(v) f)(x)
&=({\cal A}^\dagger f)(x) - \divg(\sigma v f)(x).
\end{align*}

For a density $\rho$ and a function $g$, define $\langle \rho, g
\rangle := \int_{\Re^d} g(x) \rho(x) \ud x$. With this notation, define the controlled
Lagrangian ${\cal L}:C^2(\Re^d;\Re^{+})\times C^1(\Re^d;\Re^p) \times \Re \to \Re$ as follows:
\[
{\cal L}(\rho,v\,;y) := \half \langle \rho,  |\,v\;|^2 + h^2\rangle + y\,\langle \rho,\tilde{\cal A}(v)h \rangle.
\]
The justification of this form of the
Lagrangian starting from the relative entropy cost appears in
Appendix~\ref{apdx:lagrangian-sde}.


For a given fixed observation path $z = \{z_t:0\leq t\leq T\}$, the
optimal control problem is as follows:
\begin{subequations}\label{eq:opt-cont-sde}
	\begin{align}
	&\mathop{\text{Min}}_{\pi_0, u} : \sJ(\pi_0, u\,;z) ={\sf D}(\pi_0\|\nu_0) - z_T \langle \pi_T,h\rangle + \int_0^T {\cal
		L}(\pi_t,u_t;z_t) \ud t \label{eq:cost-sde}\\
	&\text{Subj.} :\frac{\partial\pi_t}{\partial t}(x) = (\tilde{\cal A}^\dagger(u_t)\pi_t)(x) .
	\end{align}
\end{subequations}

\medskip

\begin{remark}
	This optimal control problem is a mean-field-type problem on account
	of the presence of the entropy term ${\sf D}(\pi_0\|\nu_0)$ in the
	objective function.  The Lagrangian is in a standard
	stochastic control form and the problem can be solved as a stochastic
	control problem as well~\cite{mitter2003}.  In this paper, the
	mean-field-type optimal control formulation is stressed as a
	straightforward way to
	derive the equations of the nonlinear smoothing.  
\end{remark}

The solution to this problem is given in the following proposition, whose proof appears in the Appendix~\ref{apdx:opt-ctrl-sde}. 
\begin{proposition}\label{thm:opt-ctrl-sde}
	Consider the optimal control problem~\eqref{eq:opt-cont-sde}.
	For this problem, the Hamilton's equations are as follows:
	\begin{subequations}\label{eq:hamiltons-eqn-sde}
		\begin{flalign}
		&\text{(forward)}&\frac{\partial \pi_t}{\partial t}(x) &= (\tilde{\cal
			A}^\dagger(u_t) \pi_t)(x) \label{eq:hamiltons-eqn-sde-a}&\\
		&\text{(backward)}&-\frac{\partial\lambda_t}{\partial t} (x) &= e^{-(\lambda_t(x) - z_th(x))}({\cal A}e^{\lambda_t(\cdot) - z_th(\cdot)})(x) - \half h^2(x) \label{eq:hamiltons-eqn-sde-b}&\\
		&\text{(boundary)}& \lambda_T(x) &= z_Th(x). \nonumber &
		\end{flalign}
	\end{subequations}
	The optimal choice of the other boundary condition is as follows:
	\[
	\pi_0(x) = \frac{1}{\normfactor} \nu_0(x) e^{\lambda_0(x)}
	\]
	where $\normfactor = \int_{\Re^d}
	\nu_0(x)e^{\lambda_0(x)}\ud x$ is the normalization factor.
	The optimal control is as follows:
	\[
	u_t(x) = \sigma^\top(x)\, \nabla (\lambda_t-z_th)(x).
	\]
\end{proposition}


\subsection{Optimal control: finite state-space}
\label{sec:finite_ss}

The modified process $\tilde{X}$ is a
Markov chain that also evolves in $\mathbb{S} =
\{e_1,e_2,\hdots,e_d\}$.  The control problem is parametrized by the
following:
\begin{enumerate}
	\item The initial distribution denoted as $\pi_0\in\Re^d$.  
	\item The state transition matrix denoted as $\tilde{A}(v)$ where
	$v\in (\Re^+)^{d\times d}$ is the control input.  After~\cite[Sec. 2.1.1.]{van2006filtering}, it is
	defined as follows:
	\[
	[\tilde{A}(v)]_{ij} = \begin{cases}
	[A]_{ij}[v]_{ij}\quad & i\neq j\\
	-\sum_{j\neq i} [\tilde{A}(v)]_{ij} & i=j
	\end{cases}
	\]
	and we set $[v]_{ij} = 1$ if $i=j$ or if $[A]_{ij} = 0$.
\end{enumerate}

To set up the optimal control problem, define a function
$C:(\Re^+)^{d\times d}\to\Re^d$ as follows
\[
[C(v)]_i = \sum_{j=1}^d [A]_{ij}[v]_{ij}(\log [v]_{ij}-1),\quad i = 1,\ldots,d .
\]
The Lagrangian for the optimal control problem is as follows:
\[
{\cal L}(\rho,v;y) := \rho^\top(C(v)+\half h^2) + y \;\rho^\top(\tilde{A}(v)h).
\]
The justification of this form of the
Lagrangian starting from the relative entropy cost appears in
Appendix~\ref{apdx:lagrangian-finite}.


For given observation path $z=\{z_t:0\le t\le T\}$, the optimal
control problem is as follows:
\begin{subequations}\label{eq:opt-cont-finite}
	\begin{align}
	\mathop{\text{Min }}_{\pi_0, u} &:\sJ(\pi_0, u\,;z) 
	= {\sf D}(\pi_0\|\nu_0) - z_T \pi_T^\top h + \int_0^T {\cal L}(\pi_t,u_t;z_t) \ud t \label{eq:cost-finite}\\
	\text{Subj.} &:\frac{\ud \pi_t}{\ud t} = \tilde{A}^\top(u_t) \pi_t. \label{eq:dyn_finite}
	\end{align}
\end{subequations}

The solution to this problem is given in the following proposition, whose proof appears in the Appendix.

\begin{proposition}\label{thm:opt-ctrl-finite}
	Consider the optimal control problem~\eqref{eq:opt-cont-finite}. For this problem, the Hamilton's equations are as follows:
	\begin{subequations}\label{eq:hamiltons-eqn-finite}
		\begin{flalign}
		&\text{(forward)}&\frac{\ud \pi_t}{\ud t} &= \tilde{A}^\top(u_t) \pi_t \label{eq:hamiltons-eqn-finite-a}&\\
		&\text{(backward)}&-\frac{\ud \lambda_t}{\ud t} &=
		\dv(e^{-(\lambda_t-z_th)}) \; A\, e^{\lambda_t-z_th} - \half h^2 \label{eq:hamiltons-eqn-finite-b}&\\
		&\text{(boundary)}& \lambda_T &= z_Th. \nonumber &
		\end{flalign}
	\end{subequations}
	The optimal boundary condition for $\pi_0$ is given by:
	\[
	[\pi_0]_i = \frac{1}{\normfactor} [\nu_0]_i[e^{\lambda_0}]_i,\quad i=1,\ldots,d
	\]
	where $\normfactor = \nu_0^\top e^{\lambda_0}$. The optimal control is
	\[
	[u_t]_{ij} = e^{([\lambda_t - z_th]_j-[\lambda_t - z_th]_i)}.
	\]

\end{proposition}

\medskip

\subsection{Derivation of the smoothing equations}


The pathwise equations of nonlinear filtering are
obtained through a coordinate transformation.  The proof for the following proposition is contained in the Appendix \ref{apdx:forward-backward}.  

\begin{proposition}\label{thm:forward-backward}
	Suppose $(\pi_t(x),\lambda_t(x))$ is the solution to the Hamilton's equation \eqref{eq:hamiltons-eqn-sde}.  Consider the following transformation:
	\[
	\mu_t(x) = \log(\pi_t(x)) -\lambda_t(x) + \log(\normfactor).
	\]   
	The pair $(\mu_t(x), \lambda_t(x))$ satisfy path-wise smoothing equations~\eqref{eq:pw-Zakai-a}-\eqref{eq:pw-Zakai-b}.
	Also,
	\[
	{\sf P}(X_t\in\ud x\;|\clZ_T) = {\pi}_t(x)\ud x \quad\forall t\in[0,T].
	\]
	For the finite state-space case~\eqref{eq:hamiltons-eqn-finite}, the analogous formulae are as follows:
	\[
	[\mu_t]_i = \log( [\pi_t]_i) -[\lambda_t]_i +  \log(\normfactor)
	\]
	and 
	\[
	{\sf P}(X_t = e_i\;|\clZ_T) = [{\pi}_t]_i \quad\forall t\in[0,T]
	\]
	for $i = 1,\ldots,d$.
\end{proposition}



\subsection{Relationship to the log transformation}

In this paper, we have stressed the density control viewpoint.
Alternatively, one can express the problem as a stochastic control
problem for the $\tilde{X}$ process.  For this purpose, define the
cost function $l:\Re^d \times \Re^p \times \Re \rightarrow \Re $ as follows:
\[
l(x,v\,;y) := \half |v|^2 + h^2(x) + y(\tilde{\cal A}(v)h)(x).
\]
The stochastic optimal control problem for the Euclidean case then is
as follows:
\begin{subequations}\label{eq:opt-cont-sde-hjb}
	\begin{align}
	\mathop{\text{Min }}_{\pi_0, U_t} &:\sJ(\pi_0,U_t\,;z) = \E\Big(\log \frac{\ud \pi_0}{\ud \nu_0}(\tilde{X}_0) - z_T h(\tilde{X}_T) + \int_0^T l(\tilde{X}_t,U_t\,;z_t)\ud t\Big)\\
	\text{Subj.} &:\ud \tilde{X}_t = a(\tilde{X}_t)\ud t + \sigma(\tilde{X}_t)(U_t\ud t + \ud \tilde{B}_t).
	\end{align}
\end{subequations}
Its solution is given in the following proposition whose proof appears
in the Appendix~\ref{apdx:proof-hjb}.  

\begin{proposition}\label{thm:opt-ctrl-sde-hjb}
	Consider the optimal control problem~\eqref{eq:opt-cont-sde-hjb}.
	For this problem, the HJB equation for the value function $V$ is
	as follows:
	\begin{align*}
	-\frac{\partial V_t}{\partial t}(x) &= \big({\cal A}(V_t+z_th)\big)(x) + h^2(x) -\half|\sigma^\top\nabla (V_t+z_th)(x)|^2\\
	V_T(x) &= - z_Th(x).
	\end{align*}	
	The optimal control is of the state feedback form as follows:
	\[
	U_t = u_t(\tilde{X}_t) 
	\]
	where $u_t(x) = -\sigma^\top \nabla(V_t + z_th)(x)$.  
\end{proposition}

\medskip

The HJB equation thus is exactly the Hamilton's
equation~\eqref{eq:hamiltons-eqn-sde-b} and 
\[
V_t(x) = -\lambda_t(x) ,\quad \forall x\in\Re^d,\ \forall\, t\in[0,T].
\]
Noting $\lambda_t(x) = \log q_t(x) + h(x)z_t$, the HJB equation for the
value function $V_t(x)$ is related to the backward Zakai equation for
$q_t(x)$ through the log transformation
(see also~\cite[Eqn. 1.4]{FlemingMitter82}):
\[
V_t(x) = -\log \big(q_t(x)e^{z_th(x)}\big).
\]

\subsection{Linear Gaussian case}

The linear-Gaussian case is a special case in the Euclidean setting
with the following assumptions on the model:
\begin{enumerate}
	\item The drift is linear in $x$.  That is, 
	\[ 
	a(x)=A^\top x \;\; \text{and}\;\; h(x) = H^\top x
	\] 
	where $A\in\Re^{d\times d}$ and $H\in\Re^{d}$.
	\item The coefficient of the process noise 
	\[
	\sigma(x) = \sigma
	\] 
	is a constant matrix.  We denote $Q:=\sigma\sigma^\top\in\Re^{d\times
		d}$.
	\item The prior $\nu_0$ is a Gaussian distribution with mean
	$\bar{m}_0\in\Re^{d}$ and variance $\Sigma_0\succ 0$.
\end{enumerate}

For this problem, we make the following restriction: The control input
$u_t(x)$ is restricted to be constant over $\Re^d$.  That is, the
control input is allowed to depend only upon time.  With such a
restriction, the controlled state evolves according to the sde:
\[
\ud \tilde{X}_t = A^\top \tilde{X}_t \ud t + \sigma u_t \ud t + \sigma \ud \tilde{B}_t,\quad \tilde{X}_0\sim {\cal N}(m_0,V_0).
\]  
With a Gaussian prior, the distribution $\pi_t$ is also Gaussian whose
mean $m_t$ and variance $V_t$ evolve as follow:
\begin{align*}
\frac{\ud m_t}{\ud t} &= A^\top m_t + \sigma u_t\\
\frac{\ud V_t}{\ud t} &= A^\top V_t + V_t A + \sigma\sigma^\top.
\end{align*}
Since the variance is not affected by control, the only constraint for
the optimal control problem is due to the equation for the mean.  

It is an easy calculation to see that for the linear model,
\[
(\tilde{\cal A}(v) h)(x) = H^\top (A^\top x + \sigma v).
\]
Therefore, the Lagrangian becomes 
\begin{align*}
{\cal L}(\rho,v;y) & = |v|^2 + |H^\top m|^2 + \tr(HH^\top V) + y H^\top (A^\top m + \sigma v)
\end{align*}
provided that $\rho \sim {\cal N}(m,V)$. 

For Gaussian distributions $\pi_0 = {\cal N}(m_0,V_0)$ and $\nu_0={\cal N}(\bar{m}_0,\Sigma_0)$, the divergence is given by the well known formula
\[
	{\sf D}(\pi_0\|\nu_0) = \half\log\frac{|V_0|}{|\Sigma_0|} - \frac{d}{2} + \half \tr(V_0\Sigma_0^{-1}) + \half(m_0-\bar{m}_0)^\top \Sigma_0^{-1}(m_0-\bar{m}_0)
\]
and the term due to the terminal condition is easily evaluated as
\[
\langle \pi_T, h\rangle = H^\top m_T.
\]
Because the control input does not affect the variance process,
we retain only the terms with mean and the control and express the
optimal control problem as follows: 
\begin{subequations}\label{eq:opt-cont-linear}
	\begin{align}
	&\mathop{\text{Minimize}}_{m_0, u}: \sJ(m_0,u\,;z) =  \half(m_0-\bar{m}_0)^\top {\Sigma}_0^{-1}(m_0-\bar{m}_0) \label{eq:opt-cont-linear-a} \\
	&\quad\quad\quad+ \int_0^T\half |u_t|^2 + \half|H^\top m_t|^2 + z_t^\top H^\top \dot{m}_t \ud t-z_T^\top H^\top m_T \nonumber \\
	&\text{Subject to} : \frac{\ud m_t}{\ud t} = A^\top m_t + \sigma u_t. \label{eq:opt-cont-linear-b}
	\end{align}
\end{subequations}
By a formal integration by parts, 
\begin{align*}
\sJ(m_0,u\,;z) &= \half(m_0-\bar{m}_0)^\top \bar{\Sigma}_0^{-1}(m_0-\bar{m}_0) \\
&+ \int_0^T \half |u_t|^2 + \half|\dot{z} - H^\top m_t|^2\ud t - \int_0^T\half|\dot{z}_t|^2\ud t.
\end{align*}
This form appears in the construction of the minimum energy
estimator~\cite[Ch. 7.3]{bensoussan2018estimation}.

\section{Conclusions}\label{sec:conc}

In this paper, we provide a self-contained
exposition of the equations of nonlinear smoothing as well as connections and
interpretations to some of the more recent developments in
mean-field-type optimal control theory.  These connections suggest
that the numerical approaches for mean-field
type optimal control problems can also be applied to obtain approximate filters.
Development of numerical techniques, e.g.,
particle filters to empirically
approximate the conditional distribution, has been an area of intense
research interest;
cf.,~\cite{reich2019data} and references therein.  
Approximate particle
filters based upon approximation of dual optimal control-type problems
have appeared in~\cite{chetrite2015,kappen2016adaptive,reich2019data,ruiz_kappen2017,sutter2016variational}.


\section{Appendix}

\subsection{Derivation of Lagrangian: Euclidean case}\label{apdx:lagrangian-sde}

By Girsanov's theorem, the Radon-Nikodym derivative is obtained
(see~\cite[Eqn. 35]{reich2019data}) as follows:
\[
\frac{\ud \tilde{P}}{\ud P}(\tilde{X}) = \frac{\ud \pi_0}{\ud \nu_0}(\tilde{X}_0) \; \exp\Big(\int_0^T \half|u_t(\tilde{X}_t)|^2 \ud t + u_t(\tilde{X}_t) \ud \tilde{B}_t \Big).
\]
Thus, we obtain the relative entropy formula:
\begin{align*}
{\sf D}(\tilde{P}\|P) &= \E\Big(\log\dfrac{\ud \pi_0}{\ud \nu_0}(\tilde{X}_0) + \int_0^T \half|u_t(\tilde{X}_t)|^2 \ud t + u_t(\tilde{X}_t) \ud \tilde{B}_t \Big)\\
&={\sf D}(\pi_0\|\nu_0) + \int_0^T \half\langle \pi_t,|u_t|^2\rangle \ud t.
\end{align*}

\subsection{Derivation of Lagrangian: finite state-space  case}\label{apdx:lagrangian-finite}

The derivation of the Lagrangian is entirely analogous to the Euclidean case except the R-N derivative is given according to~\cite[Prop. 2.1.1]{van2006filtering}:
\begin{align*}
\frac{\ud \tilde{P}}{\ud P}(\tilde{X}) &= \frac{\ud \pi_0}{\ud \nu_0}(\tilde{X}_0)\exp\Big(-\sum_{i,j} \int_0^T [A]_{ij}[u_t]_{ij}1_{\tilde{X}_t = e_i} \Big)\\
&\quad\quad\quad \prod_{0<t\leq T} \sum_{i\neq j}[u_{t-}]_{ij} 1_{\tilde{X}_{t-} = e_i}1_{\tilde{X}_{t} = e_j}.
\end{align*}
Upon taking log and expectation of both sides,  
we arrive at the relative entropy formula:
\begin{align*}
{\sf D}(\tilde{P}\|P) &= \E\Big(\log\dfrac{\ud \pi_0}{\ud \nu_0}(\tilde{X}_0) + \int_0^T -\sum_{i,j}[A]_{ij}[u]_{ij} 1_{\tilde{X}_t = e_i}\Big)\\
& \quad + \E\Big(\sum_{0<t\leq T} \sum_{i\neq j} \log [u_{t-}]_{ij} 1_{\tilde{X}_{t-} = e_i}1_{\tilde{X}_{t} = e_j}\Big)\\
&={\sf D}(\pi_0\|\nu_0) + \int_0^T \pi_t^\top C(u_t) \ud t.
\end{align*}

\subsection{Proof of Proposition~\ref{thm:opt-ctrl-sde}}
\label{apdx:opt-ctrl-sde}

The standard approach is to incorporate the constraint into the objective function by introducing the Lagrange multiplier $\lambda = \{\lambda_t:0\leq t\leq T\}$ as follows:
\begin{align*}
\tilde{J}(u,\lambda\,;\pi_0,z)&= {\sf D}(\pi_0 \| \nu_0) + \int_0^T \half \langle \pi_t, |u_t|^2 + h^2\rangle + z_t\langle \pi_t, \tilde{\cal A}(u_t)h\rangle \ud t\\
&\quad +\int_0^T \langle \lambda_t,  \frac{\partial\pi_t}{\partial t} - \tilde{\cal A}^\dagger(u_t) \pi_t \rangle   \ud t - z_T\langle \pi_T,h\rangle.
\end{align*}
Upon using integration by parts and the definition of the adjoint operator, after some manipulation involving completion of squares, we arrive at
\begin{align*}
\tilde{\sf J}(u,&\lambda\,;\pi_0,z)={\sf D}(\pi_0 \| \nu_0) + \int_0^T \half \langle \pi_t, |u_t - \sigma^\top\nabla(\lambda_t - z_th)|^2\rangle \ud t\\
&-\int_0^T \langle \pi_t, \frac{\partial}{\partial t}\lambda_t + {\cal A}(\lambda_t - z_th)-\half h^2+\half|\sigma^\top\nabla(\lambda_t - z_th)|^2\rangle \ud t\\
&+ \langle \pi_T,\lambda_T-z_Th\rangle - \langle \pi_0,\lambda_0\rangle.
\end{align*}
Therefore, it is natural to pick $\lambda$ to satisfy the following partial differential equation:
\begin{align}
-\frac{\partial\lambda_t}{\partial t}(x) &= \big({\cal A} (\lambda_t(\cdot) - z_th(\cdot))\big) - \half h^2(x)+\half \big|\sigma^\top\nabla (\lambda_t - z_th)(x)\big|^2 \label{eq:lambda-basic-form}\\
&= e^{-(\lambda_t(x) - z_th(x))}({\cal A}e^{\lambda_t(\cdot) - z_th(\cdot)})(x) - \half h^2(x) \nonumber
\end{align}
with the boundary condition $\lambda_T(x) = z_Th(x)$. With this choice, the objective function becomes
\begin{align*}
\tilde{\sf J}(u\,;\lambda,\pi_0,z) = {\sf D}(\pi_0 \| \nu_0) - \langle \pi_0, \lambda_0\rangle + \int_0^T \half\pi_t\big( \big|u_t -  \sigma^\top \nabla(\lambda_t-z_th)\big|^2\big) \ud t 
\end{align*}
which suggest the optimal choice of control is:
\begin{equation*}
u_t(x) = \sigma^\top(x) \nabla(\lambda_t-z_th)(x).
\end{equation*}
With this choice, the objective function becomes 
\begin{align*}
{\sf D}(\pi_0\|\nu_0) - \langle \pi_0, \lambda_0\rangle &= \int_{\mathbb{S}}\pi_0(x) \log\frac{\pi_0(x)}{\nu_0(x)}- \lambda_0(x)\pi_0(x)\ud x\\
&=\int_{\mathbb{S}} \pi_0(x)\log\frac{\pi_0(x)}{\nu_0\exp(\lambda_0(x))}\ud x
\end{align*}
which is minimized by choosing 
\[
\pi_0(x) = \frac{1}{\normfactor} \nu_0(x)\exp(\lambda_0(x))
\]
where $C$ is the normalization constant.

\subsection{Proof of Proposition~\ref{thm:opt-ctrl-finite}}
\label{apdx:opt-ctrl-finite}

The proof for the finite state-space case is entirely analogous to the proof for the Euclidean case. The Lagrange multiplier $\lambda = \{\lambda_t\in\Re^d: 0\leq t\leq T\}$ is introduced to transform the optimization problem into an unconstrained problem:
\begin{align*}
\tilde{\sf J}(u,\lambda\,;\pi_0,z) &= {\sf D}(\pi_0\|\nu_0)+ \int_0^T\pi_t^\top \big(C(u_t)+\half h^2 + z_t \tilde{A}(u_t)h \big)\ud t \\
& + \int_0^T\lambda_t^\top\big(\frac{\ud \pi_t}{\ud t} - \tilde{A}^\top(u_t) \pi_t\big)\ud t - z_T h^\top \pi_T.
\end{align*}
Upon using integral by parts,
\begin{align*}
\tilde{\sf J}(u,\lambda\,;\pi_0,z)&= {\sf D}(\pi_0\|\nu_0) + \int_0^T \pi_t^\top \big(C(u_t) -\tilde{A}(u_t)(\lambda_t - z_th)\big) \ud t\\
&+\int_0^T\pi_t^\top (-\dot{\lambda}_t+ \half h^2)\ud t +\pi_T^\top(\lambda_T-z_Th) - \pi_0^\top \lambda_0.
\end{align*}
The first integrand is
\begin{align*}
[C(u_t) -\tilde{A}(u_t)(\lambda_t - Z_th)]_i 
&= \sum_{j\neq i} A_{ij}\big([u]_{ij}(\log[u_t]_{ij}-1)\\
&\quad -[u_t]_{ij}([\lambda_t - Z_th]_j-[\lambda_t - Z_th]_i)\big)- A_{ii}.
\end{align*}
The minimizer is obtained, element by element, as
\[
[u_t]_{ij} = e^{([\lambda_t - z_th]_j-[\lambda_t - z_th]_i)}
\]
and the corresponding minimum value is obtained by:
\[
[C(u_t^*) -\tilde{A}_t(\lambda_t - Z_th)]_i = -[Ae^{\lambda_t-z_th}]_i [e^{-(\lambda_t-z_th)}]_i.
\]
Therefore with the minimum choice of $u_t$ above,
\begin{align*}
\tilde{\sf J}(u,\lambda\,;\pi_0,z)&= {\sf D}(\pi_0\|\nu_0) + \int_0^T \pi_t^\top \big(-(Ae^{\lambda_t-z_th}) \cdot  e^{-(\lambda_t-z_th)}\big) \ud t\\
&+\int_0^T\pi_t^\top (-\dot{\lambda}_t+ \half h^2)\ud t +\pi_T^\top(\lambda_T-z_Th) - \pi_0^\top \lambda_0.
\end{align*}
Upon choosing $\lambda$ according to:
\[
-[\dot{\lambda}_t]_i = [Ae^{\lambda_t-z_th}]_i [e^{-(\lambda_t-z_th)}]_i - \half h_i^2,\quad \lambda_T = z_Th.
\]
The objective function simplifies to 
\[
{\sf D}(\pi_0\|\nu_0) - \pi_0^\top \lambda_0 = \sum_{i=1}^d[\pi_0]_i \log\frac{[\pi_0]_i}{[\nu_0]_ie^{[\lambda_0]_i}}
\]
where the minimum value is obtained by choosing
\[
[\pi_0]_i =\frac{1}{\normfactor} [\nu_0]_ie^{[\lambda_0]_i}
\]
where $C$ is the normalization constant.  

\subsection{Proof of Proposition~\ref{thm:forward-backward}}
\label{apdx:forward-backward}

\paragraph{Euclidean case} Equation~\eqref{eq:hamiltons-eqn-sde-b} is
identical to the backward path-wise 
equation~\eqref{eq:pw-Zakai-b}.  So, we need to only derive the
equation for $\mu_t$.  
Using the regular form of the product formula,
\begin{align*}
\frac{\partial \mu_t}{\partial t} &= \frac{1}{\pi_t}\frac{\partial \pi_t}{\partial t} -\frac{\partial \lambda_t}{\partial t}\\
&=\frac{1}{\pi_t}(\tilde{\cal A}^\dagger(u_t)\pi_t) + e^{-(\lambda_t - z_th)}({\cal A}e^{\lambda_t(\cdot) - z_th(\cdot)}) - \half h^2.
\end{align*}
With optimal control $u_t = \sigma^\top\nabla(\lambda_t-z_th)$,
\begin{align*}
(\tilde{\cal A}^\dagger(u_t)\pi_t)
&=({\cal A}^\dagger\pi_t)-\divg\big(\sigma\sigma^\top \nabla
\pi_t\big)+\pi_t \divg\big(\sigma\sigma^\top \nabla(\mu_t+z_th)\big)\\
&\quad +(\nabla \pi_t)^\top(\sigma\sigma^\top \nabla(\mu_t+z_th))
\end{align*}
and
\begin{align*}
e^{-(\lambda_t - z_th)}({\cal A}e^{\lambda_t(\cdot) - z_th(\cdot)})
&=\frac{1}{\pi_t}({\cal A}\pi_t) - \half|\sigma^\top \nabla \log \pi_t|^2 - ({\cal A}(\mu_t + z_th)) \\
&\quad + \half \big|\sigma^\top\nabla \log(\pi_t) - \sigma^\top \nabla(\mu_t +z_th)\big|^2.
\end{align*}
Therefore,
\begin{align*}
\frac{\partial \mu_t}{\partial t} 
&=\frac{1}{\pi_t}\big(({\cal A}^\dagger\pi_t) + ({\cal A}\pi_t)-\divg(\sigma\sigma^\top \nabla \pi_t)\big)
-({\cal A}(\mu_t + z_th)) \\
&\quad + \divg\big(\sigma\sigma^\top \nabla(\mu_t+z_th)\big)+ \half\big|\sigma^\top \nabla(\mu_t+z_th)\big|^2 - \half h^2\\
&=e^{-(\mu_t(x)+z_th(x))}\big({\cal A}^\dagger
e^{(\mu_t(\cdot) +z_th(\cdot) )}\big)(x) - \half h^2(x)
\end{align*}
with the boundary condition $\mu_0 = \log \nu_0$.

\paragraph{Finite state-space case}
Equation~\eqref{eq:hamiltons-eqn-finite-b} is identical to the
backward path-wise equation~\eqref{eq:pw-Zakai-finite-b}.
To derive the equation for $\mu_t$, use the product formula
\begin{align*}
\Big[\frac{\ud \mu_t}{\ud t}\Big]_i &= \frac{1}{[\pi_t]_i}\Big[\frac{\ud \pi_t}{\ud t}\Big]_i - \Big[\frac{\ud \lambda_t}{\ud t}\Big]_i\\
&=\frac{1}{[\pi_t]_i}\big[\tilde{A}^\top(u_t)\pi_t\big]_i +  [e^{-(\lambda_t-z_th)}]_i[Ae^{\lambda_t+z_th}]_i - \half [h^2]_i.
\end{align*}
The first term is:
\begin{align*}
\big[\tilde{A}^\top(u_t)\pi_t\big]_i 
&= \sum_{j=1}^d \Big([A]_{ji}[u_t]_{ji}[\pi_t]_j-[A]_{ij}[u_t]_{ij}[\pi_t]_i\Big)
\end{align*}
and the second term is:
\begin{align*}
[e^{-(\lambda_t-z_th)}]_i&[Ae^{\lambda_t+z_th}]_i \\
&= \frac{1}{[\pi_t]_i}[e^{\mu_t+z_th}]_i\sum_{j=1}^d [A]_{ij} [\pi_t]_j [e^{-(\mu_t +z_th)}]_j.
\end{align*}
The formula for the optimal control gives
\begin{align*}
[u_t]_{ij} 
&=\frac{[\pi_t]_j}{[\pi_t]_i}[e^{-(\mu_t + z_th)}]_j[e^{\mu_t + z_th}]_i.
\end{align*}
Combining these expressions,
\begin{align*}
\Big[\frac{\ud \mu_t}{\ud t}\Big]_i &= \sum_{j=1}^d[A]_{ji}[e^{-(\mu_t + z_th)}]_i[e^{\mu_t + z_th}]_j - \half [h^2]_i\\
&=[e^{-(\mu_t + z_th)}]_i [A^\top e^{\mu_t + z_th}]_i - \half [h^2]_i
\end{align*}
which is precisely the path-wise form of the equation~\eqref{eq:pw-Zakai-finite-a}. 
At time $t=0$, $\mu_0 = \log(\normfactor[\pi_0]_i) - [\lambda_0]_i = \log[\nu_0]_i
$.

\paragraph{Smoothing distribution}
Since $(\lambda_t, \mu_t)$ is the solution to the path-wise form of
the Zakai equations, the optimal trajectory
\[
\pi_t = \frac{1}{\normfactor}e^{\mu_t+\lambda_t}
\]
represents the smoothing distribution.

\subsection{Proof of Proposition~\ref{thm:opt-ctrl-sde-hjb}}\label{apdx:proof-hjb}

The dynamic programming equation for the optimal control problem is given by
(see~\cite[Ch. 11.2]{bensoussan2018estimation}):
\begin{equation}\label{eq:HJB-V}
\min_{u\in\Re^p} \Big\{\frac{\partial V_t}{\partial t}(x) + (\tilde{\cal A}(u) V_t)(x) + l(x,u\,;z_t)\Big\} = 0.
\end{equation}
Therefore,
\begin{align*}
-\frac{\partial V_t}{\partial t}(x) &= ({\cal A}V_t)(x) + h^2(x) + z_t({\cal A}h)(x) \\
&+ \min_{u}\Big\{\half |u|^2 + u^\top\big(\sigma^\top \nabla V_t(x) + z_t\sigma^\top \nabla h(x)\big)\Big\}.
\end{align*}
Upon using the completion-of-square trick, the minimum is attained by a feedback form:
\[
u^* = -\sigma^\top \nabla (V_t +z_th)(x).
\]
The resulting HJB equation is given by
\begin{align*}
-\frac{\partial V_t}{\partial t}(x) &= \big({\cal A}(V_t + z_th)\big)(x) + h^2(x) -\half|\sigma^\top\nabla (V_t+z_th)|^2
\end{align*}
with boundary condition $V_T(x) = - z_Th(x)$. 
Compare the HJB equation with the equation~\eqref{eq:lambda-basic-form} for $\lambda$, and it follows
\[
V_t(x) = -\lambda_t(x).
\]

\bibliographystyle{spmpsci}
\bibliography{duality,tempbib}

\end{document}